\documentclass[a4paper]{article}
\usepackage[latin1]{inputenc}
\usepackage{amscd}
\usepackage[english]{babel}
\usepackage{amsmath,amssymb,amsfonts}
\usepackage{syntonly}

\begin{document}

\title{Combinatorial Morse theory and  minimality of hyperplane arrangements}
\author {M. Salvetti \thanks{ Dipartimento di
  Matematica "L. Tonelli", Largo B. Pontecorvo 5, 56127 Pisa, Italy
  (partially supported by M.U.R.S.T. 40\%)}
  \and S. Settepanella \thanks {Dipartimento di Matematica "L. Tonelli", Largo B.
  Pontecorvo 5, 56127, Pisa, Italy}}
\date {January 2007}
\maketitle

\newcommand{\F}[0]{\mathbb{F}}
\newcommand{\Z}[0]{\mathbb{Z}}
\newcommand{\R}[0]{\mathbb{R}}
\newcommand{\C}[0]{\mathbb{C}}
\newcommand{\Q}[0]{\mathbb{Q}}
\newcommand{\N}[0]{\mathbb{N}}
\newcommand{\K}[0]{\mathbb{K}}
\newcommand{\St}[0]{\mathcal S}
\newcommand{\Tbf}[0]{\mathbf{T}}
\newcommand{\Cal}[1]{\mathcal{#1}}
\newcommand{\vs}[0]{\vspace{3mm}}
\newcommand{\ph}[0]{\varphi}
\newcommand{\la}[0]{\lambda}
\newcommand{\bd}[1]{\textbf{#1}}
\newcommand{\e}{\bd e}
\renewcommand{\S}{\bd S}
\renewcommand{\lessdot}{\vartriangleleft}
\newcommand{\MS}{\tilde{\bd S}}
\newcommand{\te}{\theta}
\newcommand{\Te}{\mathbf{\theta}}
\newcommand{\q}[1]{\mbox{\bfseries{\textit{#1}}}}
\newcommand{\itl}[1]{\textit{#1}}
\newcommand{\g}[1]{\mathfrak{#1}}
\newcommand{\p}[0]{p}
\newcommand{\geets}[0]{\longleftarrow}
\newcommand{\too}[0]{\longrightarrow}
\newcommand{\sst}[0]{\subset}
\newcommand{\cl}[1]{\mathcal{#1}}
\newcommand{\G}[0]{\mathcal{G}}
\newcommand{\into}[0]{\hookrightarrow}
\newcommand{\codim}[0]{\mbox{codim }}
\newcommand{\diag}[0]{\Sigma}
\newcommand{\ea}[0]{\underline{\epsilon}}
\newcommand{\eb}[0]{\overline{\epsilon}}
%\mbox{\footnotesize{diag} }}
\newcommand{\im}[0]{\mbox{im }}
\newcommand{\az}[0]{\"{a}}
\newcommand{\ls}[0]{\mathcal L}
\newcommand{\bin}[2]{  \left(  \begin{array}{c}  #1 \\ #2  \end{array}
  \right)  }
\newcommand{\qbin}[2]{  \left[  \begin{array}{c}  #1 \\ #2
    \end{array}  \right]  }
\newcommand{\qed}[0]{\hspace{\stretch{1}}$\Box$}
\newcommand{\eq}[1][r]
       {\ar@<-3pt>@{-}[#1]
        \ar@<-1pt>@{}[#1]|<{}="gauche"
        \ar@<+0pt>@{}[#1]|-{}="milieu"
        \ar@<+1pt>@{}[#1]|>{}="droite"
        \ar@/^2pt/@{-}"gauche";"milieu"
        \ar@/_2pt/@{-}"milieu";"droite"}
\newcommand{\imm}[1][r] {\ar@{^{(}->}[#1]}
\newcommand{\D}[0]{D^{\scriptscriptstyle{(0)}}}
\newcommand{\dd}[0]{d^{\scriptscriptstyle{(0)}}}
\renewcommand{\ni}[0]{\noindent}
\newcommand{\Ga}[0]{\Gamma}
\newtheorem{df}{Definition}[section]
\newtheorem{teo}{Theorem}
\newtheorem{prop}[df]{Proposition}
\newtheorem{lem}[df]{Lemma}
\newtheorem{cor}[df]{Corollary}
\newtheorem{rmk}[df]{Remark}
\newtheorem{notat}[df]{Notation}
\newtheorem{ex}{Example}
\newtheorem{pf}{Proof}
\newtheorem{for}{}
\newenvironment{es}[1][Example.]{\begin{trivlist}
     \item[\hskip \labelsep {\bfseries #1}]}{\end{trivlist}}
\newenvironment{dm}[1][Proof.]{\begin{trivlist} \item[\hskip
    \labelsep {\bfseries #1}]}{\end{trivlist}}
\newenvironment{os}[1][Remark.]{\begin{trivlist}
     \item[\hskip \labelsep {\bfseries #1}]}{\end{trivlist}}
\newenvironment{grafi}

\section{Introduction}
In \cite{dimcapapa}, \cite{randell} it was proven that the
complement to a hyperplane arrangement in  $\C^n$ is a {\it minimal}
space, i.e. it has the homotopy type of a $CW$-complex with exactly
as many $i$-cells as the $i$-th Betti number $b_i.$ The arguments
use (relative) Morse theory and  Lefschetz type theorems.

This result of "existence" was refined in the case of complexified
real arrangements in \cite{yoshinaga}. The author consider a  flag
$V_0\subset V_1\subset\cdots\subset V_n\subset \R^n,$ $dim(V_i)=i,$
which is {\it generic} with respect to the arrangement, i.e. $V_i$
intersects transversally all codimensional$-i$ intersections of
hyperplanes. The interesting main result is a correspondence between
the $k-$cells of the minimal complex and the set of {\it chambers}
which intersect $V_k$ but do not intersect $V_{k-1}.$ The arguments
still use the Morse theoretic proof of Lefschetz theorem, and some
analysis of the critical cells is given. Unfortunately, the
description does not allow to understand exactly the attaching maps
of the cells of a minimal complex.

In this paper we give, for a complexified real arrangement $\mathcal
A$, an  explicit Lefschetz theorem - free description of a minimal
$CW$-complex. The idea is that, since an explicit $CW-$complex $\S$
which describes the homotopy type of the complement already exists
(see \cite{salvetti}), even if not minimal, one can work over such
complex trying to "minimize" it. A natural tool for doing that is to
use {\it combinatorial Morse theory} over $\S.$ We follow the
approach of \cite{forman}, \cite{forman1} to combinatorial Morse
theory (i.e., Morse theory over $CW$-complexes).

So, we explicitly construct  a combinatorial gradient vector field
over $\S,$ related to a given system of polar coordinates in $\R^n$
which is {\it generic} with respect to the arrangement $\mathcal A.$
Let $\St$ be the set of all {\it facets} of the stratification of
$\R^n$ induced by the arrangement $\mathcal A$ (see
\cite{bourbaki}). Then $\St$ has a natural partial ordering given by
$F\prec G$ iff $clos(F)\supset G.$ Our definition of genericity of a
coordinate system, which is stronger than that used in
\cite{yoshinaga}, allows to give a {\it  total ordering} $\lessdot$
on $\St,$ which we call the {\it polar ordering} of the facets.

 The $k$-cells in $\S$ are in one-to-one correspondence with
the pairs $[C\prec F^k]$, where $C$ is a chamber in $\St$ and $F^k$
is a codimensional-$k$ facet of $\St$ which is contained in the
closure of $C.$ Then the gradient field can be recursively defined
as the set of pairs
$$([C\prec F^{k-1}],[C\prec F^k])$$
such that $F^{k-1}\prec F^k$ and $F^k\lessdot F^{k-1},$ and such
that the origin cell of the pair is not the end cell of another pair
of the field. We also give a non-recursive equivalent
characterization of the field (thm. \ref{T6}, (ii)) only in terms of
the partial ordering $\prec$ and of the total ordering $\lessdot.$

Analog to index$-k$ critical points  in the standard Morse theory,
there are  {\it singular cells} of dimension $k$: they are those
$k$-cells which do not belong to the gradient field (see
\cite{forman}). In our situation, they are given (see corollary
\ref{cor:sing1}) in terms of the orderings by those $[C\prec F^k]$
such that
\medskip

\ni $\begin{array}{lcccccc}
       i) & F^k \lessdot F^{k+1}, & \forall & F^{k+1} & s.t. & F^k\prec F^{k+1}\ \\
       ii) & F^{k-1}  \lessdot  F^k, & \forall &  F^{k-1} & s.t. &
       C\prec F^{k-1}\prec F^k \
     \end{array}$

\medskip

It is easy to see that associating to a singular cell $[C\prec F^k]$
the unique chamber $C'$ which is {\it opposite} to $C$ with respect
to $F^k,$ gives a one-to-one correspondence between the set of
singular cells in $\S$ and the set of all chambers in $\St.$ So we
also derive by our method the main result in \cite{yoshinaga}.

Minimality property of the complement follows easily, so the above
singular cells of $\S$ give an explicit basis for the integral
cohomology, which depends on the system of polar coordinates (we
call such a basis a {\it polar basis} for the cohomology). A minimal
complex is obtained from $\S$ by contracting all pairs of cells
which belong to the vector field.

Our construction gives also an explicit algebraic complex which
computes local system cohomology of $M(\mathcal A).$ In dimension
$k$ such complex has one generator for each singular cell of $\S.$
The boundary operator is obtained by a method which is the
combinatorial analog to "integrating over all paths" which satisfy
some conditions. We give a reduced formula for the boundary, which
is effectively computable in terms only of the two orderings
$\prec,\ \lessdot.$ For abelian local systems, the boundary operator
assumes an even  nicer reduced form. There exists a vast literature
about calculation of local system cohomology on the complement to an
arrangement: several people constructed algebraic complexes
computing local coefficient cohomology, in the abelian case (see for
example \cite{dcohen, cohenorlik, esnschvie, kohno, libyuz,
schtervar, suciu, yoshinaga}). Our method seems to be more effective
than the previous ones.

In the last part we find a generic polar ordering on the braid
arrangement. We give a description of the complex $\S$ in this case
in terms of {\it tableaux} of a special kind; next, we characterize
the {\it singular} tableaux and we find an algorithm to compare two
tableaux with respect to the polar ordering.

Some of the most immediate remaining problems are: first, compare
polar bases with the well-known $nbc-$bases of the cohomology (see
\cite{bjorner_ziegler, orlik_terao}); second, characterize polar
orderings in a purely combinatorial way (so, using an {\it oriented
matroid} counterpart of generic polar coordinates). \footnote{A
preliminary version of this paper was published as a preprint in
\cite{salvettisette}}

\section{n-dimensional polar coordinates}\label{sec1}

 For reader's convenience, we recall here $n$-dimensional polar
coordinates.  Since usually one knows only standard 3-dimensional
formulas, we give here coordinate changes in general.

Start with an orthonormal basis
$$\e_1,...,\e_n$$
of the Euclidean $n-$dimensional space $V$ and let
$$P\ \equiv \ (x_1,...,x_n)$$
the associated cartesian coordinates of a point $P.$\ We will
confuse the point $P$ and the vector $OP,$ \ $O$ being the origin of
the coordinate system.

Let in general
$$pr_W:\ V\ \to\ W$$
be the orthogonal projection onto a subspace $W$ of $V.$ Consider
the two flags of subspaces
$$V_i\ =\ <\e_1,...,\e_i>, \quad i=0,...,n\ (V_0={0})$$
and
$$W_i\ =\ <\e_i,...,\e_n>, \quad i=1,...,n.$$
Let
$$P_i\ :=\ pr_{W_i}(P),\ i=1,...,n$$
(so $P_1=P$). One has
$$P_i\ =\ pr_{W_i}(P_j), \ j\leq i$$
so there are orthogonal decompositions

$$ P_{i}\ =\ P_{i+1}\ +\ x_i \e_i,\ x_i\in\R,\ i=1,...,n \eqno(1)$$
(set $P_{n+1}=0$)

 Clearly
$$\begin{array}{rcl}
P_i\ = \ 0 & \Rightarrow & P_j=0\ \text{for}\ j\geq i \\
P_i\ \neq 0 & \Rightarrow & P_j \neq 0\ \text{for}\ j\leq i
\end{array}\eqno(2)$$

Let
$$\te_{n-1}\in (-\pi,\pi]$$
be the angle that $OP_{n-1}$ forms with $e_{n-1}$ (in the 2-plane
$W_{n-1}$). Let then
$$\te_i\in [0,\pi],\quad i=1,...,n-2$$
be the angle that $OP_i$ makes with $e_i.$

The polar coordinates of $P$ will be given by the "module"
$$\rho\ =\ \|P\|$$
more (if $P\neq 0$) "arguments"
$$\te_{1},...,\te_{n-1}$$
(defined only for $i\leq max\{j:\ P_j\neq 0\}$).

The coordinate change between polar and cartesian coordinates is
given by

$$\begin{array}{rcl}
x_1 & = & \rho\cos(\te_1)\\
x_2 & = & \rho\sin(\te_1)\cos(\te_2)\\
% x_3 & = & \rho\sin(\te_1)\sin(\te_2)\cos(\te_3)\\
\vdots & \vdots & \vdots \\
x_i & = & \rho\sin(\te_1)...\sin(\te_{i-1})\cos(\te_i)\\
\vdots & \vdots & \vdots \\
% x_{n-2} & = & \rho \sin(\te_1)...\sin(\te_{n-3})\cos(\te_{n-2})\\
x_{n-1} & = & \rho \sin(\te_1)...\sin(\te_{n-2})\cos(\te_{n-1}) \\
x_n & = & \rho \sin(\te_1)...\sin(\te_{n-1})\\

\end{array} \eqno(3)$$
Notice that these formulas make sense always if we conventionally
set $\te_i= 0$ for $P_i=0.$

The inverse formulas are

$$\begin{array}{rcl}
\rho^2 & = & x_1^2+...+x_n^2 \\
\cos^2(\te_{1}) & = & \frac{x_{1}^2}{x_1^2+...+x_{n}^2 }\\
% \cos^2(\te_{2}) & = & \frac{x_{2}^2}{x_2^2+...+x_{n}^2 }\\
\vdots & \vdots & \vdots \\
\cos^2(\te_{i}) & = & \frac{x_{i}^2}{x_i^2+...+x_{n}^2 }\\
\vdots & \vdots & \vdots \\
% \cos^2(\te_{n-2}) & = & \frac{x_{n-2}^2}{x_{n-2}^2+...+x_{n}^2 }\\
\cos^2(\te_{n-1}) & = & \frac{x_{n-1}^2}{x_{n-1}^2+x_{n}^2 }\\
\end{array} \eqno(4)$$

\section{Combinatorial Morse theory}\label{CM}

We recall here the main points of Morse theory for $CW$-complexes,
from a combinatorial viewpoint. All the definitions and results in
this section are taken from \cite{forman}, \cite{forman1}.

We restrict to the case of our interest, that of {\it regular}
$CW-$complexes.

\bigskip

\subsection{Discrete Morse functions}
Let $M$ be a finite regular $CW$-complex, let $K$ denote the set of
cells of $M,$ partially ordered by
$$\sigma\ <\ \tau \quad \Leftrightarrow\quad \sigma\subset \tau,$$
and $K_p$ the cells of dimension p.

\begin{df} A discrete Morse function on $M$ is a function
$$f:K \too \R$$
satisfying for all $\sigma^{(p)} \in K_p$ the two conditions

$$\begin{array}{cccccc} (i) && \sharp \{\tau^{(p+1)}>\sigma^{(p)}\ |\ f(\tau^{(p+1)})
\leq f(\sigma^{(p)})\} & \leq & 1& \\
(ii) &&  \sharp \{v^{(p-1)} < \sigma^{(p)} \ |\ f(\sigma^{(p)})
\leq f(v^{(p-1)}) \} & \leq & 1&\\
\end{array}$$

We say that $\sigma^{(p)} \in K_p$ is a $\mathrm {critical\ cell}$
of index $p$ if the cardinality of both these sets is $0$.
\end{df}

\begin{rmk} One can show that, for any given cell of $M,$ at least one of the two
cardinalities in (i),\ (ii) is $0$ (\cite{forman}).
\end{rmk}

Let $m_p(f)$ denote the number of critical cells of $f$ of index
$p$. As in the standard theory one can show (see \cite{forman})

\begin{prop} M is homotopy equivalent to a $CW$-complex with exactly $m_p(f)$ cells of dimension p.
\end{prop}

\subsection{Gradient vector fields}

Let $f$ be a discrete Morse function on a $CW$-complex $M$. One can
define the discrete gradient vector field $V_f$ of $f$ as:
$$V_f\ =\ \{(\sigma^{(p)},\tau^{(p+1)}) | f(\tau^{(p+1)}) \leq f(\sigma^{(p)})\}.$$
By definition of Morse function, each cell belongs to at most one
pair of $V_f.$ More generally, one defines

\begin{df} A discrete vector field\ $V$\ on \ $M$\ is a collection of pairs
$(\sigma^{(p)},\tau^{(p+1)})$ of cells such that each cell of $M$
belongs to at most one pair of $V$.
\end{df}

Given a discrete vector field $V$ on $M$, a $V$-path is a sequence
of cells
\begin{equation}\label{sequenza}\sigma_0^{(p)},\tau_0^{(p+1)},
\sigma_1^{(p)},\tau_1^{(p+1)},\sigma_2^{(p)},\cdots,\tau_r^{(p+1)},\sigma_{r+1}^{(p)}
\end{equation}
such that for each $i=0,\cdots ,r$ $(\sigma_r^{(p)},\tau_r^{(p+1)})
\in V$ and
$\sigma_i^{(p)} \neq \sigma_{(i+1)}^p < \tau_i^{(p+1)}$.\\
Such a path is a non trivial closed path if $0 \leq r$ and
$\sigma_0^{(p)} = \sigma_{(r+1)}^p$.  One has:

\begin{teo}\label{teo:gradfield} A discrete vector field V is the gradient vector field
of a discrete Morse function if and only if there are no non-trivial
closed V-path.
\end{teo}

\begin{rmk}
An equivalent combinatorial definition of discrete vector field is
that of {$\mathrm {matching}$} over the Hasse diagram of the poset
associated to the CW-complex (see for ex. \cite{forman1}).
\end{rmk}

\bigskip

\section{Applications to Hyperplane arrangements}

\subsection{Notations and recalls}\label{S-complex} \ Let $\Cal A\ =\ \{H\}$
be a finite affine hyperplane
arrangement in $\R^n.$  Assume $\Cal A$ essential, so that the
minimal dimensional non-empty intersections of hyperplanes are
points (which we call \emph{vertices} of the arrangement).
Equivalently, the maximal elements  of the associated
\emph{intersection lattice} $L(\Cal A)$ (see \cite{orlik_terao})
have rank $n.$

Let
$$M(\Cal A)\ =\ \C^n \setminus \bigcup_{H\in A}\ H_{\C}$$
be the complement to the complexified arrangement. We use the
regular CW-complex $\S\ =\ \S(\Cal A)$ constructed in
\cite{salvetti} which is a deformation retract of $M(\Cal A)$ (see
also \cite{gelfand_rybnikov},\ \cite{bjorner_ziegler},\
\cite{orlik_terao},\ \cite{salvetti2}). Here we recall very briefly
some notations and properties.

Let
$$\Cal S:=\{F^k\}$$
be the stratification of $\R^n$ into facets $F^k$ which is induced
by the arrangement (see \cite{bourbaki}), where exponent $k$ stands
for codimension.  Then $\Cal S$ has standard partial ordering
$$F^i \ \prec F^j \quad \text{iff}\quad clos(F^i)\supset F^j$$
Recall that  $k$-cells of $\S$ bijectively correspond to pairs
$$[C\prec F^k]$$
where $C=F^0$ is a chamber of $\Cal S.$

Let $|F|$ be the affine subspace spanned by $F,$ and let us consider
the subarrangement
$$\Cal A_F\ =\ \{ H \in \Cal A\ :\ F\subset H\}.$$
A cell $[C\prec F^k]$ is in the boundary of $[D\prec G^j]$ ($k<
j$) iff

\bigskip

i) $F^k\prec G^j$

ii) the chambers $C$ and $D$ are contained in the same chamber of
$\Cal A_{F^k}.$

\bigskip

 Previous conditions are equivalent to say that
  $C$ is the chamber of $\Cal A$ which is "closest" to $D$
among those which contain $F^k$ in their closure.
\begin{notat}\label{nt:boundary} i) We denote the chamber $D$ which
appear in the boundary cell $[D\prec G^{j}]$
of a cell $[C\prec F^k]$  by  $C.G^{j}.$

ii) More generally, given a chamber $C$ and a facet $F,$ we denote
by $C.F$ the unique chamber containing $F$ and lying in the same
chamber as $C$ in $\Cal A_{F^k}.$ Given two facets $F,\ G$ we will
use also for  $(C.F).G$ the notation (without brackets) $C.F.G$.
\end{notat}

It is possible to realize $\S$ inside $\C^n$ with explicitly given
attaching maps of the cells (see \cite{salvetti}).  Recall also that
the construction can be given for any \emph{oriented matroid} (see
the above cited references).

\subsection{Generic polar coordinates}\label{sec:polarcoord}
In general,  we distinguish between \emph{bounded} and
\emph{unbounded} facets. Let $B(\Cal S)$ be the union of bounded
facets in $\Cal S.$ When $\Cal A$ is central and essential (i.e
$\cap_{H\in\Cal A}\ H$ is a single point $O\in V$) then $B(\Cal
S)=\{O\}.$ In general, it is known that $B(\Cal S)$ is a compact
connected subset of $V$  and the closure of a small open
neighborhood $U$ of $B(\Cal S)$ is homeomorphic to a ball (so $U$ is
an open ball; see for ex. \cite{salvetti}).

Given a system of polar coordinates associated to $O,\e_1,...,\e_n,$
the coordinate subspace $V_i\ , \ i=1,...,n $ (see section 1) is
divided by $V_{i-1}$ into two components:
$$V_i\setminus V_{i-1}\ =\ V_i(0)\cup V_i(\pi)$$
where
$$V_i(0)\ =\ \{P:\ \te_{i}(P)=0\}$$
and
$$V_i(\pi)\ =\ \{P:\ \te_i(P)=\pi\}$$
More generally, we indicate by
$$V_i(\bar{\te}_i,...,\bar{\te}_{n-1})\ :=\ \{P:\ \te_i(P)=\bar{\te}_i,...,\te_{n-1}(P)=\bar{\te}_{n-1}\}
\eqno(5)$$
where by convention $\bar{\te}_j=0\ \text{or}\ \pi\ \Rightarrow\
\bar{\te}_k=0$ for all $k>j;$ \ so in particular,
$V_i(0)=V_i(0,...,0)$ and $V_i(\pi)=V_i(\pi,0,...,0)$ ($n-i$
components). The space $V_i(\bar{\te})$ is an $i-$dimensional open
half-subspace in the euclidean space $V,$ and we denote by
$|V_i(\bar{\te})|$ the subspace which is spanned by it. We have from
(3)
$$|V_i(\bar{\te}_i,...,\bar{\te}_{n-1})| \ =\
<\e_1,...,\e_{i-1},\bar{\e}> $$ where
$$\bar{\e}=\bar{\e}(\te_i,...,\te_{n-1}):=\ \sum_{j=i}^n\
\left(\prod_{k=i}^{j-1}\ \sin(\te_k)\right)\ \cos(\te_{j})\ \e_j$$

 For all $\delta\in (0,\pi/2)$ the space
$$\tilde{B}:=\ \tilde{B}(\delta):=\ \{P:\ \te_i(P)\in(0,\delta),\
i=1,...,n-1,\ \rho(P)>0 \}$$ is an open cone contained in $\R^n_+.$

\begin{df}\label{df:gencoor}
We say that a system of polar coordinates in $\R^n,$  defined by an
origin $O$ and a base $\e_1,...,\e_n,$ is \emph{generic} with
respect to the arrangement $\Cal A$ if it satisfies the following
conditions:
\bigskip

i) the origin $O$ is contained in a chamber $C_0$ of $\Cal A;$

ii) there exist $\delta\in (0,\pi/2)$  such that
$$B(\Cal S)\ \subset\ \tilde{B}=\tilde{B}(\delta);$$
(therefore, for each facet $F\in\St$ one has
$F\cap\tilde{B}\neq\emptyset$);

iii) subspaces
$V_i(\bar{\te})=V_i(\bar{\te}_{i},...,\bar{\te}_{n-1})$ which
intersect $clos(\tilde{B})$ (so $\bar{\te}_j\in[0,\delta]$ for
$j=i,...,n-1$) are generic with respect to $\Cal A,$ in the sense
that, for each codim$-k$ subspace $L\in L(\Cal A),$
$$i\geq k \ \Rightarrow\ V_i(\bar{\te})\cap L\cap clos(\tilde{B})\neq\emptyset\ \text{and}\
dim (|V_i(\bar{\te})|\cap L)\ =\ i-k.$$
\end{df}
\bigskip

It is easy to see that genericity condition implies that the origin
$O$ of coordinates belongs to an unbounded chamber. It turns out
that such chamber must intersects the infinity hyperplane
$H_{\infty}$ into a relatively open set. This is equivalent to say
that the sub-arrangement given by the walls of the chamber is
essential.

In fact we have

\begin{teo} For each unbounded chamber $C$ such that $C\cap H_{\infty}$
is relatively open, the set of points $O\in C$ such that there
exists a polar coordinate system centered in $O$ and generic with
respect to $\Cal A$ forms an open subset of $C.$
\end{teo}
\bigskip
\ni {\bf Proof.} We proceed by first proving  the following
\medskip
\begin{lem} Let $\Cal A$ be a central essential arrangement in $V.$ Then
there exist orthonormal frames $\e_1,...,\e_n$ which are
\emph{generic} with respect to $\Cal A,$ in the sense that each
subspace $V_i:=<\e_1,...,\e_i>,\ i=1,...,n,$ intersects
transversally each $L\in L(\Cal A).$ Given a chamber $C,$ the first
vector $\e_1$ can be any vector inside $C.$

Actually, the set of generic frames is open inside the space of
orthonormal  frames in $V.$
\end{lem}
\medskip
\ni {\bf Proof of lemma.} Let $O'$ be the intersection of all
hyperplanes, and take an orthonormal coordinate system with basis
$\e'_1,...,\e'_n.$ Then each hyperplane $H$ is given by a linear
form
$$H_i\ =\ \{x: (\alpha_i\ \cdot \ x)\ =\ 0\},\ i=1,...,|\Cal
A|$$ where we denote by $(\ \cdot\ )$ the canonical inner product.
Any codimensional$-k$ $L$ is given by an intersection of $k$
linearly independent hyperplanes $H_{i_1},...,H_{i_k}$ of $\Cal A.$
The genericity condition on a frame $\e_1,...,\e_n$ is written as
$$rk\left[(\alpha_{i_r}\ \cdot\ \e_s)\right]_{{\scriptstyle{
r=1,...,k}}\atop \scriptstyle{s=1,...,i}}\ =\ min\{k,i\}$$ or,
equivalently
$$rk\left[(\alpha_{i_r}\ \cdot\ \e_s)\right]_{{\scriptstyle{
r=1,...,k}}\atop \scriptstyle{s=1,...,k}}\ =\ k.\ \eqno(7)$$ It is
clear that genericity applied to $V_1$ gives that $\e_1$ is not
contained in any hyperplane, i.e. it belongs to some chamber of
$\Cal A.$
 Equation (7) is easily translated into the
equivalent one
$$dim(<\alpha_{i_1},...,\alpha_{i_k},\e_{k+1},...,\e_n>)=n. \eqno(8)$$

Passing to the dual space $V^*$ by using the inner product, the set
of all hyperplanes
$$<\alpha_{i_1},...,\alpha_{i_{n-1}}>\subset V^*$$
gives an arrangement $\Cal A^*.$ Since $\e_1$ belongs to a chamber
of $\Cal A,$ each subspace \ $<\alpha_{i_1},...,\alpha_{i_k}>$
intersects transversally the orthogonal $\e_1^{\bot},$ so $\Cal A^*$
induces an arrangement $\Cal A_1:=\Cal A^*\cap \e_1^\bot$ over
$\e_1^\bot.$ Condition (8) requires an orthonormal basis
$\e_2,...,\e_n\in \e_1^\bot$ which is generic with respect to the
flag
$$V'_1:=<\e_n>,...,V'_i:=<\e_{n-i+1},...,\e_n>,... $$
Then we conclude the proof of the first and second assertions  by
induction on $n.$

For the last one, notice that
$$\e_1,\ \e_n,\ \e_2,\ \e_{n-1},\ ...$$
can vary respectively in a chamber of the arrangements
$$\Cal A,\ \Cal A_1=\Cal A^*\cap \e_1^\bot,\ \Cal A_{1,n}:=(\Cal
A_1)^*\cap\e_n^\bot,\ \Cal A_{1,n,2}:=\Cal A_{1,n}^*\cap \e_2^\bot,\
...$$ which is an open set inside orthonormal frames. \qed
\bigskip

We come back to the prove of theorem.
\medskip

\ni \emph{Case I}:\ $\Cal A$ central and essential.

\ni Let $O'$ be the center of $\Cal A.$ According to the previous
lemma we can find $\e_1,...,\e_n$ generic with respect to $\Cal A,$
and with $\e_1:=\frac{OO'}{\|OO'\|}.$ If we consider a system of
polar coordinates associated to $O,\e_1,...,\e_n$ then the subspaces
$V_i$ satisfy condition (iii) of genericity. Perform a small
translation onto $\Cal A,$
$$x_i\ \rightarrow \ x_i+\sigma,\ 0<\sigma<<1$$
which moves the center $O'$ into the positive octant. Then if
$$\sigma<<\delta<<1$$
all conditions in definition 3.1 are satisfied by continuity and the
fact that genericity is an open condition.

\bigskip

\ni \emph{General case.}

\ni In case of an affine arrangement referred to a system of
cartesian coordinates $O',\e'_1,...,\e'_n,$ hyperplanes are written
as
$$H_i=\{x:\ (\alpha_i\ \cdot\ x)\ =\ a_i\}\ ,i=1,...,|\Cal A|.$$
Let
$$H_i^0:=\{(\alpha_i\ \cdot\ x)\ =\ 0)\}$$
be the \emph{direction} of $H_i$ and let $\Cal A_0$ be the
associated central arrangement. Notice that if $\Cal A$ is essential
then so is $\Cal A_0.$ We can assume without loss of generality that
$\|\alpha_i\|=1,\ \forall i,$ so the vector
$$a_i\cdot\alpha_i$$
represents the translation taking $H_i^0$ into $H_i.$

Let $C$ be an unbounded chamber of $\Cal A$ such that $C\cap
H_{\infty}$ is relatively open in $H_{\infty}.$ Then the directions
of the walls of $C$ are the walls of a chamber $C'$ in $\Cal A_0.$
By previous case, there exist points $O$ in  $C'$  and systems of
polar coordinates $O, \e_1,...,\e_n$ which are generic with respect
to $\Cal A_0.$ Let $\delta>0$ satisfy definition \ref{df:gencoor}
for one of such systems. We can assume (up to a homotethy of center
$O'$)
$$|a_i|<<\delta, \ \forall i.$$
Then the same system satisfies the definition for $\Cal A.$

Of course, the condition of genericity is open, so we finish the
proof of the theorem.  \qed

\subsection{Orderings on $\St$}\label{sec:orderings}
Fix a system of generic polar coordinates, associated to a center
$O$ and frame $\e_1,...,\e_n.$ \ Let $\delta>0$ be the number coming
from definition 3.1. We denote for brevity \
$\bar{B}:=clos(\tilde{B}(\delta)).$\  Each point $P$ has polar
coordinates $P\ \equiv\ (\te_0,\te_1,...,\te_{n-1}),$ where we use
the convention $\te_0:=\rho.$

We remark that when the pole $O$ is very far, the polar coordinates
of one point inside $\tilde{B}(\delta)$ are approximately the same
as its standard Cartesian coordinates.

Notice that (5) makes sense also for $i=0,$ being
  $$V_0(\bar{\te}_0,\bar{\te}_1,...,\bar{\te}_{n-1})$$
given by a single point $P$ with
$$\rho(P)=\bar{\te}_0,\ \te_1(P)=\bar{\te}_1,\ ...\ ,\ \te_{n-1}(P)=\bar{\te}_{n-1}.$$

Given a codimensional$-k$ facet $F\in \St,$ let us denote by
$$F(\te):=F(\te_i,...,\te_{n-1}):=F\cap V_i(\te_i,...,\te_{n-1}),\quad
\te_j\in[0,\delta],\ j=i,...,n-1 $$ (notice: \ $F=F(\te)=F\cap V_n$\
with $\te=\emptyset$.)

 By genericity conditions, if \ $i\geq k$\
then \ $F(\te)$ \ is either empty or it is a codimensional $k+n-i$
facet contained in $V_i(\te_i,...,\te_{n-1}).$

Let us set, for every facet $F(\te),$
$$i_{F(\te)}\ :=\ min\{j\geq 0:\ V_{j}\cap clos(F(\te))\neq\emptyset\}.$$
Still by genericity, setting $L:=|F(\te)|,$ one has
$$L\cap V_{j}\ \neq\ \emptyset\ \Leftrightarrow\ j\geq codim(F(\te))$$
so also
$$i_{F(\te)}\ \geq \ codim(F(\te)) \eqno(9).$$

When the facet $F(\te):=F(\te_{i},...,\te_{n-1}),\ i>0,$\ is not
empty and $i_{F(\te)}\geq i$ \ (i.e., $clos(F(\te))\cap
V_{i-1}=\emptyset$), then among its vertices ($0-$dimensional facets
in its boundary) there exists, still by genericity,  a unique one
$$P:=P_{F(\te)}\in clos(F(\te)) \eqno(10)$$
such that
$$\te_{i-1}(P)\ =\ min\{\te_{i-1}(Q):\ Q\in
clos(F(\te))\}\eqno(11)$$ (of course, $P_{F(\te)}=F(\te)$ if
$dim(F(\te))=0,$\ i.e. $i=k$).

When $i_{F(\te)}<i$ then the point $P$ of (10) is either the origin
$0$ ( $\Leftrightarrow \ i_{F(\te)}=0$ $\Leftrightarrow$ $F$ is the
base chamber $C_0$) or it is the unique one such that
$$\te_{i_{F(\te)}-1}(P)\ =\ min\{\te_{i_{F(\te)}-1}(Q):\ Q\in
clos(F(\te))\cap V_{i_{F(\te)}}\} \eqno(12)$$

\bigskip

\begin{df}\label{df:minpoint}
Given any facet $F(\te)=F(\te_i,...,\te_{n-1})$ let us denote by
$$P_{F(\te)}\in clos(F(\te))$$
the "minimum" vertex of \ $clos(F(\te))\cap\ V_{i_{F(\te)}}$ \ (as
in (10))

\ni (for $F\in \St$ we briefly write \ $P_F$).

\medskip

We associate to the facet $F(\te)$ the $n-$vector of polar
coordinates of $P_{F(\te)}$
$$\Theta(F(\te))\ :=\ (\te_0(F(\te)),...,\te_{i_{F(\te)}-1}(F(\te)),0,...,0)$$
($n-i_{F(\te)}$ zeroes) where we set
$$\te_j(F(\te)):=\te_j(P_{F(\te)}),\quad j=0,...,i_F-1.$$
\end{df}

Notice that when $i_{F(\te)}\geq i$ then all coordinates $\te_j$ of
$\te$ with $j\geq i_{F(\te)}$ must be zero.

We want to define another ordering over the poset $\St,\prec.$ We
give a recursive definition, actually ordering all facets in
$V_i(\te)$ for any given $\te=(\te_1,...,\te_{n-1}).$

\bigskip

\begin{df}[Polar Ordering]\label{df:maindef} Given $F,\ G\ \in\ \St,$ and
given $\bar{\te}=(\bar{\te}_i,...,\bar{\te}_{n-1}),$ $0\leq i\leq
n,\ \bar{\te}_j\in [0,\delta]$ for $j\in i,..,n-1,$
($\bar{\te}=\emptyset$ for $i=n$) such that $F(\bar{\te}),
G(\bar{\te}) \neq \emptyset$, we set
$$F(\bar{\te})\ \lessdot\ G(\bar{\te})$$
iff one of the following cases occur:

\medskip

i)\ $P_{F(\bar{\te})}\neq P_{G(\bar{\te})}.$\ Then\
$\Theta(F(\bar{\te}))\ <\ \Theta(G(\bar{\te}))$ \ according to the
anti-lexicographic ordering of the coordinates (i.e., the
lexicographic ordering starting from the last coordinate).

\bigskip

ii)\ $P_{F(\bar{\te})}=P_{G(\bar{\te})}.$\ Then either

\medskip
\ni iia) $dim(F(\bar{\te}))=0$ \ (so
$P_{F(\bar{\te})}=F(\bar{\te})$) and $F(\bar{\te})\neq G(\bar{\te})$
(so $dim(G(\bar{\te}))>0$)

\medskip
or

\medskip
\ni iib) $dim(F(\bar{\te}))>0,\ dim(G(\bar{\te}))>0.$ \ In this case
let \ $i_0:=i_{F(\bar{\te})}=i_{G(\bar{\te})}.$

 When $i_0\geq i$ (case (11)) one can write
$$\Theta(F(\bar{\te}))=\Theta(G(\bar{\te}))=(\tilde{\te}_0,...,\tilde{\te}_{i-1},\bar{\te}_i,...,
\bar{\te}_{i_0-1},0,...,0).$$
 Then \ $\forall\epsilon,\ 0<\epsilon<<\delta,$\ it must happen
$$F(\tilde{\te}_{i-1}+\epsilon,\bar{\te}_i,...,\bar{\te}_{i_0-1},0,...,0)\ \lessdot\
G(\tilde{\te}_{i-1}+\epsilon,\bar{\te}_i,...,\bar{\te}_{i_0-1},0,...,0).$$

If $i_0<i$ (as in (12)) then one can write
$$\Theta(F(\bar{\te}))=\Theta(G(\bar{\te}))=(\tilde{\te}_0,...,\tilde{\te}_{i_0-1},0,...,0).$$
 Then \ $\forall\epsilon,\ 0<\epsilon<<\delta,$\ it must happen
$$F(\tilde{\te}_{i_0-1}+\epsilon,0,...,0)\ \lessdot\ G(\tilde{\te}_{i_0-1}+\epsilon,0,...,0).$$

($n-i_0$ zeroes) \qed

\end{df}

\bigskip

Condition (iib) says that one has to move a little bit the suitable
$V_j(\te')$ which intersects $clos(F(\te))$ and $clos(G(\te))$ in a
point $P(F(\te))=P(G(\te))$ (according to (11) or (12)), and
consider the facets which are obtained by intersection with this
"moved" subspace.
\medskip

It is quit clear from the definition that irriflessivity and
transitivity hold for $\lessdot$ so we have

\bigskip

\begin{teo}   Polar ordering $\lessdot$ is a total ordering on the
facets of $V_i(\bar{\te}),$ for any given
$\bar{\te}=(\bar{\te_i},...,\bar{\te}_{n-1}).$
 In particular
(taking $\bar{\te}=\emptyset$) it gives a total ordering on $\St.$
\qed \end{teo}
\bigskip

\bigskip

The following property, comparing polar ordering with the partial
ordering $\prec,$ will be very useful.

\begin{teo}\label{thm:locord} Each codimensional-$k$ facet $F^k\in \St$ ($k<n$) such that $F^k\cap V_k =\emptyset$
has the following property: among all codimensional-$(k+1)$ facets
$G^{k+1}$ with $F^k\prec G^{k+1},$ there exists a unique one
$F^{k+1}$ such that
$$F^{k+1}\lessdot F^k.$$
If $F^k\cap V_k\neq\emptyset$ (so $F^k\cap V_k=P(F^k)$) then
$$F^k\lessdot G^{k+1},\quad \forall G^{k+1}\ with\ F^k\prec
G^{k+1}.$$
\end{teo}
\bigskip

\ni {\bf Proof.}  In the latter case, where $F^k\cap V_k = P(F^k),$
for every facet $G^{k+1}$ in the closure of $F^k$ one has
$P(G^{k+1})\not\in clos(V_k)$ (by (9)), so $F^k\lessdot G^{k+1}.$

In general, for all facets $G$ contained in the closure of $F^k,$
one has either $P(G)\neq P(F^k)$ and $\Theta(F^k)<\Theta(G),$ so
$F^k\lessdot G,$ or $P(G)=P(F^k). $ For those $G^{k+1}$ such that
$P(G^{k+1})=P(F^k)$ one reduces, after $\epsilon-$deforming (may be
several times) like in definition \ref{df:maindef}, to the case
where $F$ is a one-dimensional facet contained in some $V_h\setminus
V_{h-1},$ with $h\geq 1,$ and for such case the assertion is clear.
\qed

\bigskip

Let $\St^{(k)}:= \St\ \cap\ V_k$ be the stratification induced onto
the coordinate subspace $V_k.$ A codimensional-$j$ facet in $V_k$ is
the intersection with $V_k$ of a unique codimensional-$j$ facet in
$\St$,\ $j\leq k.$ Let $\lessdot_k$ be the polar ordering of
$\St^{(k)},$ induced by the polar coordinates associated to the
basis $\e_1,...,\e_k$ of $V_k.$ By construction, for all $F,\
G\in\St$ which intersect $V_k,$ one has
$$F\cap V_k\ \lessdot_k\ G\cap V_k\quad \text{iff}\quad F\lessdot
G.$$ So we can say that $\lessdot_k$ is the {\it restriction} of
$\lessdot$ to $V_k$ and also $\lessdot_k$ is the restriction of
$\lessdot_h$ for $k<h.$

By genericity conditions,  for each $F^k\in\St$ there exists a
unique  $F^k_0$ with the same support and intersecting $V_k$ (in one
point).

The following recursive characterization of the polar ordering will
be used later. The proof is a direct consequence of definition
\ref{df:maindef} and theorem \ref{thm:locord}.

\begin{teo}\label{thm:algo} Assume that, for all $k=0,...,n,$ we know
the polar ordering of all the $0$-facets (=codimensional-$k$ facets)
of $\St^{(k)}$ (in particular, $\forall\ F^k$ we know whether
$F^k\cap V_k\neq\emptyset$). Then we can reconstruct the polar
ordering of all $\St.$ Assuming we know it for all facets of
codimension $\geq k+1,$ then given $F^k, \ G^k$ we have:

- if both $F^k,\ G^k$ intersect $V_k$ then the ordering is the same
as the restriction to $\St^{(k)};$

- if one intersects $V_k$ and the other does not, the former is the
lower one;

- if no of the two facets intersects $V_k,$ then let $F'^{(k+1)}$,\
(resp.\ $G'^{(k+1)}$) be the facet in the boundary of $F^k$ (resp.
$G^k$) which is minimum with respect to $\lessdot.$ Then
$$F^k \lessdot G^k$$
iff either
$$F'^{(k+1)} \lessdot G'^{(k+1)}$$
 or
 $$F'^{(k+1)} = G'^{(k+1)}\quad \text{but}\quad G_0^{(k)} \lessdot
 F_0^{(k)}$$
 where $F_0^{(k)}$ (resp. $G_0^{(k)}$) means (as above) the unique facet with the same
 support which intersects $V_k.$

 Moreover, each $F^k$ intersecting $V_k$ is lower than any
 codimensional-$(k+1)$ facet. If $F^k$ does not intersect $V_k$ then
 $F^k$ is bigger than its {\it minimal boundary} $F'^{(k+1)}$ and
 lower than any codimensional -$(k+1)$ facet which is bigger than $F'^{(k+1)}.$

 This determines the polar ordering of all facets of codimension
 $\geq k.$
\end{teo}\ \qed

\bigskip

\begin{rmk} We ask whether it is possible to characterize polar
orderings in {\it purely combinatorial} ways. The problem is more or
less that of finding a good "combinatorial" description of a flag of
subspaces (or better, of half-subspaces)  which corresponds to a
generic system of polar coordinates, so that we are able to decide
what facets belong to coordinate half-spaces. It seems quite
reasonable that this can be done by specifying combinatorially a
"generic" flag in the given oriented matroid.
\end{rmk}

\subsection{Combinatorial vector fields}
We consider here the regular CW-complex $\S=\S(\Cal A)$ of section
3.1. Recall that $k-$cells correspond to pairs $[C\prec F^k],$ where
$C$ is a chamber and $F^k$ is a codimensional-k facet in $\St.$ We
will define a combinatorial gradient vector field \ $\Phi$ over
$\S.$ One can describe $\Phi$ (see section 2.2) as a collection of
pairs of cells
$$\Phi\ =\ \{(e,f)\in \S \times \S\ |\ dim(f)=dim(e)+1,\ e\in \partial(f)\}$$
so that $\Phi$ decomposes into its dimensional-p components
$$\Phi\ =\ \bigsqcup_{p=1}^n \ \Phi_p,\qquad \Phi_p\subset \S_{p-1}\times \S_{p}$$
($\S_p$ being the $p-$skeleton of $\S$). Let us indicate by
$$\ea\ ,\ \eb:\ \Phi\ \to\ \S,    \quad \ea(a,b)=a,\ \eb(a,b)=b$$
the first and last cells of the pairs of $\Phi.$

We give the following recursive definition:

\bigskip

\begin{df}[Polar Gradient]{\label{df:maindef1}} We define a
combinatorial gradient field $\Phi$ over $\S$ in the following way:
\medskip

\ni the $(j+1)-$th component $\Phi_{j+1}$ of $\Phi,$ $j=0,...,n-1,$
is given by all pairs
$$([C\prec F^j],[C\prec F^{j+1}]),\quad F^j\prec F^{j+1}$$
(same chamber $C$) such that
\medskip

1.\ $F^{j+1}\lessdot F^{j}$

\medskip

2.\ $\forall F^{j-1}\prec F^j$ the pair
$$([C\prec F^{j-1}],[C\prec F^{j}]) \not\in \ \Phi_j$$

\end{df}

\bigskip

Notice that condition (ii) in 3.1 is automatically verified for
pairs as in definition \ref{df:maindef1}. Condition 2 of
\ref{df:maindef1} is empty for the 1-dimensional part $\Phi_1$ of
$\Phi,$ so
$$\Phi_1=\{([C\prec C],[C\prec F^{1}])\ :\quad  F^1\lessdot C\}.$$

\ni According to the definition of generic polar coordinates, only
the base-chamber $C_0$ intersects the origin $O=V_0,$ so by theorem
\ref{thm:locord} all $0-$cells $[C\prec C],\ C\neq C_0,$ belong to
exactly one pair of $\Phi_1.$

\bigskip

\begin{teo}\label{T6} One has:

\ni (i) $\Phi$ is a combinatorial vector field on $\S$ which is the
gradient of a combinatorial Morse function (according to part 2.2).
\smallskip

\ni (ii) The pair
$$([C\prec F^j],[C\prec F^{j+1}]),\quad F^j\prec F^{j+1}$$
belongs to $\Phi$ iff the following conditions hold:

  (a) $F^{j+1}\lessdot F^j$

  (b) $\forall\ F^{j-1}$ such that $C\prec F^{j-1}\prec F^j,$\ one has $F^{j-1}\lessdot F^j.$
\smallskip

\ni (iii) Given $F^{j}\in \St,$ there exists a chamber $C$ such that
the cell  $[C\prec F^{j}]\in \bar{\epsilon}(\Phi)$ iff there exists
$F^{j-1}\prec F^{j}$ with $F^{j}\lessdot F^{j-1}.$ More precisely,
for each chamber $C$ such that there exists $F^{j-1}$ with
$$C\prec F^{j-1} \prec F^{j},\quad F^{j}\lessdot F^{j-1} \eqno(*)$$
the pair \ $([C\prec \bar{F}^{j-1}],[C\prec F^{j}])\in \Phi,$ \
where $\bar{F}^{j-1}$ is the maximum $(j-1)-$facet (with respect to
polar ordering) satisfying conditions (*).
\smallskip

\ni (iv) The set of $k-$dimensional singular cells is given by
$$Sing_k(\S)\ =\ \{[C\prec F^k]\ :\ F^k\cap V_k \neq \emptyset,\
F^j\lessdot F^k, \ \forall \ C\prec F^j \precneqq F^k\}\eqno(13).$$
Equivalently, $F^k\cap V_k$ is the maximum (in polar ordering) among
all facets of $C\cap V_k.$
\end{teo}

\bigskip

\ni {\bf Proof}. \  Clearly $\Phi_1$ satisfies (ii) with $j=0.$ We
assume by induction that $\Phi_j$ is a combinatorial vector field
satisfying (ii). Consider now a $j-$cell $[C\prec F^j]\in \S.$
Assume condition (b) of (ii) holds for $F^j:$ then if there exists
$F^{j+1}$
 with $F^j\prec F^{j+1},\ F^{j+1}\lessdot F^j$ (and this happens by
 theorem \ref{thm:locord} iff $F^j\cap V_j=\emptyset$) then
 $$([C\prec F^j],[C\prec F^{j+1}]) \in \Phi_{j+1}.$$
If (b) of (ii) does not hold ($j\geq 2$) then let $F^{j-1}$ be the
biggest (according to polar ordering) codimensional $j-1$ facet such
that
$$C\prec F^{j-1}\prec F^j,\quad F^{j}\lessdot F^{j-1}.$$
Take any $F^{j-2}$ such that $C\prec F^{j-2}\prec F^{j-1}.$ We
assert that $F^{j-2}\lessdot F^{j-1}.$ Otherwise, certainly there
exists another facet  $G^{j-1}$ with
$$F^{j-2}\prec G^{j-1} \prec F^j$$
and by theorem \ref{thm:locord} it should be $F^{j-2}\lessdot
G^{j-1},$ contradicting the maximality of $F^{j-1}.$ So by induction
$$([C\prec F^{j-1}],[C\prec F^j])\ \in\ \Phi_j$$
(this proves (iii)) and the cell $[C\prec F^j]$ cannot be the origin
of a pair of $\Phi_{j+1}.$

To show that $\Phi_{j+1}$ is a vector field, we have to see that no
cell $[C\prec F^{j+1}]$ is the end of two different pairs of
$\Phi_{j+1}.$ After $\epsilon-$deforming we reduce to the case where
$F^{j+1}$ is $0-$dimensional. Then the unicity of a $j-$facet $F^j$
such that $C\prec F^j\prec F^{j+1},$ and such that (a) and (b) of
(ii) hold easily comes from convexity of the chamber $C.$

This proves both that $\Phi$ is a combinatorial vector field and
(ii).

Next, we prove that $\Phi$ is a gradient field by using theorem
\ref{teo:gradfield} of section \ref{CM}: we have to  show that
$\Phi$ has no closed loops.

So let
\bigskip

\ni $(\ [C_1\prec F^j_1],\ [C_1\prec F^{j+1}_1],\ [C_2\prec F^j_2],\
[C_2\prec F^{j+1}_2],\ ... $

 $$ ...,\ [C_m\prec F^j_m],\
[C_m\prec F^{j+1}_m],\ [C_{m+1}\prec F^j_{m+1}] \ )$$  be a
$\Phi-$path (see (1)). First,  notice that the $j+1$-facets are
ordered
$$F^{j+1}_m \ \underline{\lessdot}\ ...\ \underline{\lessdot}\ F^{j+1}_1.$$
In fact, by definition of path and the boundary in $\S$ (see sec.
\ref{S-complex} ) we have at the $k$-th step:
$$F^j_{k+1}\prec F^{j+1}_k,\quad F^j_{k+1}\prec F^{j+1}_{k+1},\quad
F^{j+1}_{k+1}\lessdot F^j_{k+1}.$$  If also
$$F^{j+1}_{k}\lessdot F^j_{k+1}$$
then by theorem \ref{thm:locord}  \ $F^{j+1}_{k+1} = F^{j+1}_{k};$ \
otherwise we have necessarily $F^{j+1}_{k+1}\lessdot
F^j_{k+1}\lessdot F^{j+1}_k.$ Then if the path is closed it derives
(still by theorem \ref{thm:locord}) that all the $F_k^{j+1}$ equal a
unique $F^{j+1}.$ Moreover, up to $\epsilon$-deforming, we can
assume that the path is contained in some $V_i(\te)$ with $F^{j+1}$
a $0-$dimensional facet. Under these assumptions, we show that
$$F_1^j\lessdot ... \lessdot F_m^j.$$

Let $V_{i-1}(\te_{i-1},\te_i,...)\subset V_i(\te)$ be the subspace
containing the point $F^{j+1};$ after $\epsilon$-deforming, the path
can be seen inside the subspace
$$\tilde{V}:=V_{i-1}(\te_{i-1}-\epsilon,\te_i,...)$$
where for each cell $[C_k\prec F^j_k]$ one has that $C_k$ is a
convex open polyhedron in $\tilde{V}$ (may be infinite) and $F^j_k$
is, by point (iii), its {\it maximum} vertex: all the facets of
$C_k$ are lower (in polar ordering) than $F^j_k.$

By the definition of boundary in \ref{S-complex} the two chambers
$C_k,\ C_{k+1}$ belong to the same chamber of $\Cal A_{F^j_{k+1}}.$
Such a chamber is a convex cone with maximum facet (with respect to
polar ordering) is $F^j_{k+1},$ and such that each of its facets has
the same support as some facet of $C_{k+1}$ of the same dimension,
having the vertex $F^j_{k+1}$ as one of its $0$-facets.  Then
clearly all the facets of $C_k$ are lower (in polar ordering) than
$F^j_{k+1}.$ In particular $F^j_k \lessdot F^j_{k+1},$ which proves
that there are no non-trivial closed $\Phi$-paths.

It remains to prove part (iv). In view of (ii), (iii), a cell
$[C\prec F^k]$ does not belong to $\Phi$ iff
$$F^k \lessdot F^{k+1},\quad \forall \ F^k\prec F^{k+1}\eqno(A)$$
and
$$F^{k-1}\lessdot F^k,\quad \forall\ C\prec F^{k-1}\prec F^k.
\eqno(B)$$ Condition (A) holds by theorem \ref{thm:locord} iff \ $P\
:=\ F^k\cap V_k\neq \emptyset.$ \ Then $P$ is a $0$-dimensional
facet in $V_k,$ and (B) holds iff $P$ is the maximum facet of the
chamber $C\cap V_k$ (according to polar ordering). This is
equivalent to (iv), and finishes the proof of the theorem. \qed

\bigskip

As an immediate corollary we have
\bigskip

\begin{cor}\label{cor:sing1} Once a polar ordering is assigned, the set of singular
cells is described only in terms of it by
\bigskip

\ni $Sing_k(\S)\ :=\ \{[C\prec F^k]\ :\ $
\medskip

\ni $\begin{array}{lcccccc}
       a) & F^k \lessdot F^{k+1}, & \forall & F^{k+1} & s.t. & F^k\prec F^{k+1} \\
       b) & F^{k-1}  \lessdot  F^k, & \forall &  F^{k-1} & s.t. &
       C\prec F^{k-1}\prec F^k \ \}
     \end{array}$

\qed
\end{cor}

\bigskip

\begin{rmk}\label{rmk:loccrit} Of course, condition b) of corollary \ref{cor:sing1}
is equivalent to:  \
$$F'\ \lessdot\ F^k \ \text{ for all}\ F' \text{ in the
interval}\ C\prec F'\ \prec\ F^k.$$
\end{rmk}

\begin{rmk}\label{rmk:sing2} By (iv) of theorem \ref{T6} $Sing_k(\S)$ corresponds to
the pairs $(C,v)$ where $C$ is a chamber of the arrangement $\Cal
A_k:=\ \Cal A\cap V_k$ and $v$ is the maximum vertex of $C.$ Then
$v$ is the minimum vertex of the chamber $\tilde{C}$ of $\Cal A_k$
which is {\it opposite} to $C$ with respect to $v.$ Of course,
$\tilde{C}\cap V_{k-1}=\emptyset,$ so we re-find the one-to-one
correspondence between the singular k-cells of $\S$ and the chambers
of $\Cal A_k$ which does not intersect $V_{k-1}$ (see
\cite{yoshinaga}).
\end{rmk}

\begin{rmk}  By easy computation, the integral boundary of the Morse
complex generated by singular cells (see \cite{forman}) is zero, so
we obtain the minimality of the complement. Alternatively, the same
result is obtained by noticing that singular cells are in one-to-one
correspondence with the set of all the chambers of $\St$ by remark
\ref{rmk:sing2}. But \ $\sum\ b_i\ =\ |\{ chambers \}|$ \ (see for
ex. \cite{zaslavsky, orlik_terao}).
\end{rmk}

\begin{rmk}
 Our description gives also an explicit additive
basis for the homology and for the cohomology in terms of the
singular cells in $\S.$ We can call it a $\mathrm{polar\ basis}$
(relative to a given system of generic polar coordinates). It would
be interesting to compare such basis with the well-known $nbc-$basis
of the cohomology (see \cite{bjorner_ziegler, orlik_terao}).
\end{rmk}

\section{Morse complex for local homology}

The gradient field indicates how to obtain a minimal complex from
$\S,$ by contracting all pairs of cells in the field. For each pair
of cells $(e^{k-1},e^k)$ in $\Phi,$ one has a contraction of $e^k$
into $\partial(e^k)\setminus int(e^{k-1}),$ by "pushing"
$int(e^{k-1})\cup int(e^k)$ onto the boundary.

In particular, it is possible to obtain a {\it Morse complex} which
computes homology and cohomology, even with local coefficients. We
describe here such an algebraic complex, computing homology with
local coefficients for the complement $M(\mathcal A).$ The boundary
operators depend only on the partial ordering $\prec$ and on the
polar ordering $\lessdot.$

First, we give to the coordinate space $V_i$ the orientation induced
by the ordered basis $\e_1,...,\e_i.$ Given a codimensional$-i$
facet $F^i\in\St,$ the support $|F^i|$ is transverse to $V_i,$ so we
give the orthogonal space $|F^i|^\perp$ the orientation induced by
that of $V_i.$ Recall from \cite{salvetti} that the complex $\S$ has
a real projection $\Re: \S \to \R^n$ which induces a
dimension-preserving cellular map onto the $dual$ cellularization
$\St^\vee\subset \R^n$ of $\St.$ We give to a cell $e(F^i)\in
\St^\vee,$ dual to $F^i,$ the orientation induced by that of
$|F^i|^\perp.$  We give to a cell $[C\prec F^i]\in \S $ the
orientation such that the real projection $\Re:[C\prec F^i]\to
e(F^i)$ is orientation preserving.

Let $L$ be a  {\it local system} over $M(\mathcal A),$ i.e. a module
over the group-algebra of the fundamental group $\pi_1(M(\mathcal
A))$  The basepoint is the origin $O\in C_0$ of the coordinates,
which can be taken as the unique $0-$cell of $\S$ (and of
$\St^\vee$) contained in $C_0.$ Up to homotopy, we can consider only
{\it combinatorial} paths in the $1-$skeleton of $\S,$ i.e.
sequences of consecutive edges. Sequences, or {\it galleries},
$$C_1,...,C_t$$
of {\it adjacent} chambers uniquely correspond to a special kind of
combinatorial paths in the $1-$skeleton of $\S,$ which we call {\it
positive} paths. Two galleries with the same ends and of minimal
length determine two homotopic positive paths (see \cite{salvetti}).
One says that a positive path, or gallery, {\it crosses} an
hyperplane $H\in \mathcal A$ if two consecutive chambers in the path
are separated by $H$.

Remark that the $1-$dimensional part $\Phi_1$ of the polar field
gives a maximal {\it tree} in the $1-$skeleton of $\S.$ Each
$0-$cell $v(C)$ of $\S$ is determined by its dual chamber $C\in
\St.$ Then each $v(C)\in \S$ is connected to the origin $O$ by a
unique path $\Gamma(C),$ which is a positive path, determined by a
gallery of chambers starting in $C$ and ending in $C_0.$ We have

\begin{lem}\label{lem:minpath} For all chambers $C,$ the path $\Gamma(C)$ is minimal,
i. e. it crosses each hyperplane at most once.
\end{lem}
\begin{dm} One has that $\Gamma(C)$ consists of a sequence of
$1-$cells $[C\prec F^1]$ where $F^1\lessdot C.$ It is sufficient to
see that the hyperplane $H=|F^1|$ separates $C$ from $C_0.$ This
comes immediately from the definition of polar ordering, since one
has $P(F^1)=P(C)$ and $F^1$ is encountered before $C$ by a half-line
$V_1(\theta_1,...).$ \qed
\end{dm}

\bigskip

\begin{notat}\label{not:path} i) Given two chambers $C,\ C'$ we denote by $\mathcal
H(C,C')$ the set of hyperplanes separating $C$ from $C'.$

ii) Given an ordered sequence of (possibly not adjacent) chambers
$C_1,...,C_t$  we denote by $u(C_1,...,C_t)$ the rel-homotopy class
of
$$u(C_1,...,C_t)\ =\ u(C_1,C_2)u(C_2,C_3)\cdots u(C_{t-1},C_t$$ where $u(C_i,C_{i+1})$
is a minimal positive path induced by a minimal gallery starting in
$C_i$ and ending in $C_{i+1}.$ We denote by
$$\bar{u}(C_1,...,C_t)\in \pi_1(M(\mathcal A),O)$$
the homotopy class of a path which is the composition
$$\bar{u}(C_1,...,C_t):=(\Gamma(C_1))^{-1}u(C_1,...,C_t)
\Gamma(C_t).$$

We denote by
$$\bar{u}(C_1,...,C_t)_*\ \in\ Aut(L)$$
the automorphism induced by $\bar{u}(C_1,...,C_t).$
\end{notat}

\bigskip

We need also some definitions.

\begin{df} A cell $[C\prec F]\in \S$ will be called {\rm locally
critical} if $F$\ is the maximum, with respect to $\lessdot,$ of all
facets in the interval $\{F'\ :\ C\prec F'\prec F\}$ of the poset
$(\St, \prec).$
\end{df}

By corollary \ref{cor:sing1} and remark \ref{rmk:loccrit} a critical
cell is also locally critical. By theorem \ref{T6}, part (iii), the
cell $[C\prec F^k]$ belongs to the $k-$dimensional part $\Phi_k$ of
the polar field iff it is not locally critical.

\begin{df} Given a codimensional-$k$ facet $F^k$ such that $F^k\cap
V_k\neq \emptyset,$ a sequence of pairwise different
codimensional-$(k-1)$ facets
$$\mathcal F(F^k) \ :=\ (F^{(k-1)}_{i_1},\cdots,F^{(k-1)}_{i_m}), \ m\geq 1$$
such that
$$F^{(k-1)}_{i_j}\prec F^{k},\ \forall \ j$$
and
$$F^k\lessdot F^{(k-1)}_{i_j}\ \text{for}\ j<m$$
while for the last element
$$F^{(k-1)}_{i_m}\lessdot F^k$$
is called an {\rm admissible $k-$sequence.}

It is called an {\rm  ordered} admissible $k-$sequence if
$$F^{(k-1)}_{i_1}\lessdot\cdots \lessdot F^{(k-1)}_{i_{m-1}}.$$
\end{df}

\bigskip

Notice that in an admissible $k-$sequence with $m=1,$ it remains
only a codimensional-$(k-1)$ facet which is lower (in polar
ordering) than the given codimensional-$k$ facet.

Two admissible $k-$sequences
$$\mathcal F(F^k) \ :=\ (F^{(k-1)}_{i_1},\cdots,F^{(k-1)}_{i_m})$$
$$\mathcal F(F'^k) \ :=\ (F'^{(k-1)}_{j_1},\cdots,F'^{(k-1)}_{j_l})$$
$F^k\neq F'^k,$ can be {\it composed} into a sequence

$$\mathcal F(F^k)\mathcal F(F'^k) \ :=\
(F^{(k-1)}_{i_1},\cdots,F^{(k-1)}_{i_m},
F'^{(k-1)}_{j_1},\cdots,F'^{(k-1)}_{j_l})$$ when for the last
element of the first one it holds
$$F^{(k-1)}_{i_m}\prec F'^k.$$
In case $F^{(k-1)}_{i_m}=F'^{(k-1)}_{j_1}$ we write this facet only
once, so there are no repetitions in the composed sequence.
\bigskip

\begin{df}\label{df:ammseq} Given a critical $k-$cell \ $[C\prec F^k]\in \S$ \ and a
critical $(k-1)-$cell \ $[D\prec G^{k-1}]\in \S,$ \ an {\rm
admissible sequence}
$$\mathcal F\ =\ \mathcal F_{([C\prec F^k],\ [D\prec G^{(k-1)}])}$$
for the given pair of critical cells is a sequence of
codimensional-$(k-1)$ facets
$$\mathcal F\ :=\ (F^{(k-1)}_{i_1},\cdots,F^{(k-1)}_{i_h})$$
obtained as composition of admissible $k-$sequences
$$\mathcal F(F^k_{j_1})\cdots \mathcal F(F^k_{j_s})$$
such that:
\medskip

a) \ $F^k_{j_1}=F^k$ \ (so $F^{k-1}_{i_1}\prec\ F^{k}$);

b) \ $F^{k-1}_{i_h}\ =\ G^{k-1}$ \ and the chamber
$$C.F^{k-1}_{i_1}.\cdots.F^{k-1}_{i_h}$$
(see notation \ref{nt:boundary}) equals \ $D$;

c) \ for all $j=1,\cdots, h$ \ the $(k-1)-$cell \
$$[C.F^{k-1}_{i_1}.\cdots.F^{k-1}_{i_j}\prec F^{k-1}_{i_j}]$$
is locally critical.

We have an {\rm ordered} admissible sequence if all the
$k-$sequences that compose it are ordered.
\end{df}

\bigskip

\begin{lem} All admissible sequences are ordered.
\end{lem}

\begin{dm} Let $s$ be an admissible sequence. One has to show that each
$k-$sequence composing $s$ is ordered. This follows by definition
\ref{df:ammseq}, c), and by the definition of polar ordering. \qed
\end{dm}

Denote by
$${\mathcal Seq}\ =\ {\mathcal Seq}([C\prec F^k],[D\prec G^{(k-1)}]$$
the set of all admissible sequences for the given pair of critical
cells. Of course, this is a finite set which is determined only by
the orderings $\prec,\ \lessdot.$ In fact, the "operation" which
associates to a chamber $C$ and a facet $F$ the chamber $C.F$ is
detected only by the Hasse diagram of the partial ordering $\prec.$
The chamber $C.F$ is  determined by: $C.F\prec F$ and $C.F$ is
connected by the shortest possible path (= sequence of adjacent
chambers) in the Hasse diagram of $\prec.$

Given an admissible sequence $s=(F^{k-1}_{i_1},..., F^{k-1}_{i_h})$
for the pair of critical cells $[C\prec F^k],$ \   $[D\prec
G^{k-1}],$ we denote (see notation \ref{not:path}) by
$$u(s)\ =\ u(C,C.F^{k-1}_{i_1},\cdots,C.F^{k-1}_{i_1}.\cdots
.F^{k-1}_{i_h}) $$ and by
$$\bar{u}(s) =\ \bar{u}(C,C.F^{k-1}_{i_1},\cdots,C.F^{k-1}_{i_1}.\cdots
.F^{k-1}_{i_h}) .$$ Set also $l(s):=h$ for the length of $s$ and
$b(s)$ for the number of $k-$sequences forming $s.$

Now we have a complex which computes local system homology.
\bigskip

\begin{teo}\label{thm:morsecompl}
The homology groups with local coefficients
$$H_k(M(\mathcal A),L)$$
are computed by the algebraic complex  \ $(C_*,\partial_*)$ such
that:

in dimension $k$
$$C_k\ :=\ \oplus\ L.\ e_{[C\prec F^k]},$$
where one has one generator for each singular cell $[C\prec F^k]$ in
$\S$ of dimension $k$.

The boundary operator is given by

$$\partial_k(l. e_{[C\prec F^k]})\ =\ \sum A^{[C\prec F^k]}_{[D\prec
G^{k-1}]}(l).\ \ e_{[D\prec G^{k-1}]} \eqno(*)$$ (\ $l\in L$\ )
where the incidence coefficient is given by:
$$A^{[C\prec F^k]}_{[D\prec G^{k-1}]}\ :=
\sum_{s\in Seq}\ (-1)^{l(s)-b(s)}\ \bar{u}(s)_*\ \eqno(**)$$ Here
the sum is over all possible  admissible sequences $s$ for the pair
$[C\prec F^k],\ [D\prec G^{k-1}].$  \qed
\end{teo}

\bigskip

\begin{dm} The proof follows by the definition of the vector field,
from theorem \ref{T6} and from the definition of boundary in $\S.$
In fact, condition c) implies (by (iii) of theorem \ref{T6}) that
the $(k-1)-$cell \ $[C.F^{k-1}_{i_1}.\cdots.F^{k-1}_{i_j}\prec
F^{k-1}_{i_j}]$ does not belong to $\Phi_{k-1},$ so the pair
$$([C.F^{k-1}_{i_1}.\cdots.F^{k-1}_{i_j}\prec
F^{k-1}_{i_j}], [C.F^{k-1}_{i_1}.\cdots.F^{k-1}_{i_j}\prec E^{k}]\
\in\ \Phi$$ for $j<h.$ The result is obtained by substituting to
$$[C.F^{k-1}_{i_1}.\cdots.F^{k-1}_{i_j}\prec F^{k-1}_{i_j}]$$ the
remaining boundary
$$\partial([C.F^{k-1}_{i_1}.\cdots.F^{k-1}_{i_j}\prec
E^{k}])\setminus [C.F^{k-1}_{i_1}.\cdots.F^{k-1}_{i_j}\prec
F^{k-1}_{i_j}]$$ and keeping into account the given orientations.
\qed
\end{dm}

\bigskip

\begin{rmk}  The sign in formula (**) can be expressed in the
following way. If $s=(F^{k-1}_{i_1},..., F^{k-1}_{i_h})$ then set
$$\alpha :=\ \#\{j<h\ : \ F^{k-1}_{i_j}\lessdot F^{k-1}_{i_{j+1}}\}$$
and set $\epsilon=0$ or $1$ according whether the first element
$F^{k-1}_{i_1}\lessdot F^k$ or $F^k\lessdot F^{k-1}_{i_1}$ Then one
has
$$(-1)^{l(s)-m(s)}\ =\ (-1)^{\alpha +\epsilon}\  .$$
\end{rmk}

\bigskip

Many admissible sequences in the boundary operator cancel, because
of the sign rule. We give a very simplified formula in the
following.

\bigskip

\begin{df} 1) Given a pair of critical cells  $[C\prec F^k],\ [D\prec G^{k-1}],$
we say that an admissible sequence
$$s=(F^{k-1}_{i_1},..., F^{k-1}_{i_h})\in Seq$$
{\rm is m-extensible by} the facet $F'^{k-1}$ if:

a)  $F'^{k-1}$ can be inserted into the sequence $s$ to form another
sequence $s'$ of length $h+1$ which is still admissible with respect
to the same pair of critical cells, and such that
$$\bar{u}(s)_*\ =\ \bar{u}(s')_*  .$$

b)  $F'^{k-1}$ is the minimum (with respect to $\lessdot$)
codimensional-$(k-1)$ facet which satisfies a) (then we call $s'$
the m-extension of $s$ by $F'^{k-1}$);

c) $F'^{k-1}$ is the minimum of the facets $F"^{k-1}$ in the
sequence $s'$ such that the sequence $s":=s'\setminus F"^{k-1}$
obtained by removing $F"^{k-1}$ is still admissible, and
$$\bar{u}(s")_*\ =\ \bar{u}(s')_* \  =\ \bar{u}(s)_*\ .$$
In other words, $s$ is not the m-extension of some $s"$ by
$F"^{k-1},$ with $F"^{k-1}\lessdot F'^{k-1}.$

2) We say that an admissible sequence
$$s=(F^{k-1}_{i_1},..., F^{k-1}_{i_h})\in Seq$$
{\rm is m-reducible by}  $F'^{k-1}$ in $s,$ if the sequence $s'$
obtained by removing $F'^{k-1}$ is m-extensible by $F'^{k-1}.$
\end{df}

\bigskip

Set \ $Seq^e$ \ and \ $Seq^r$ \ be the set of m-extensible, resp.
m-reducible (by some codimensional$-(k-1)$ facet), admissible
sequences for a given pair of critical cells. By definition
$$Seq^e\ \cap \ Seq^r\ =\emptyset.$$

The following lemma is also clear from the previous definition.

\begin{lem} There is a one-to-one correspondence
$$Seq^e\ \leftrightarrow\ Seq^r$$
which associates to a sequence $s$ which is m-extensible by
$F'^{k-1}$ its extension $s'$ (obtained by adding $F'^{k-1}$).
\end{lem}

\bigskip
Set
$$Seq^0\ :=\ Seq\setminus (Seq^e \cup  Seq^r).$$
as the set of non m-extensible and non m-reducible sequences.

Since the sign in formula (**) which is associated to an
m-extensible sequence $s$ and to its extension $s'$ is opposite, it
follows:

\begin{teo}\label{thm:morsecompl1} The coefficient of the boundary operator
in (**) of theorem \ref{thm:morsecompl} holds
$$A^{[C\prec F^k]}_{[D\prec G^{k-1}]}\ :=
\sum_{s\in Seq^0}\ (-1)^{l(s)-b(s)}\ \bar{u}(s)_*\ $$
\end{teo}\qed

\bigskip

The reduction of theorem \ref{thm:morsecompl1} is strong.

We consider now {\it abelian} local systems, i.e. modules $L$ such
that the action of $\pi_1(M(\mathcal A))$ factorizes through
$H_1(M(\mathcal A)).$ Then to each elementary loop $\gamma_H$
turning around an hyperplane $H$ in the positive sense it is
associated an element \ $t_H\in Aut(L),$ so one has homomorphisms
$$\Z[\pi_1(M(\mathcal A))]\ \to\ \Z[H_1(M(\mathcal A))]\ \to\
\Z[t_H^{\pm 1}]_{H\in\mathcal A}\subset End(L).$$ An abelian local
system as that just defined is determined by the system $\mathcal
T:=\{t_H,\ H\in \mathcal A\},$ so we denote it by $L(\mathcal T).$

Given an admissible sequence $s=(F^{k-1}_{i_1},..., F^{k-1}_{i_h})$
relative to the pair $[C\prec F^k],\ [D\prec \G^{k-1}],$ and given
an hyperplane $H\in \mathcal A,$ we indicate by $\mu(s,H)$ the
number of times the path $u(s)$ crosses $H.$

\begin{lem} For $s,\ H$ as before, one has
\medskip \

\ni 1) $ H\in\mathcal H(C_0,C)\cap\mathcal H(C_0,D)$ then
$$
\begin{array}{llll}
 \mu(s,H)=0 \ \text{if}& F^k\not\subset H \ \text{or} & F^k\subset H \ \text{and}\ F^{k-1}_{i_1}\lessdot F^k  \\
 \mu(s,H)=2 \ \text{otherwise}\\
\end{array}
$$

\ni 2) $ H\in\mathcal H(C_0,D)\cap\mathcal H(C,D)$ then $\mu(s,H)=1$

\ni 3) $ H\in\mathcal H(C_0,C)\cap\mathcal H(C,D)$ then

\ if $F^k\not\subset H$ then $\mu(s,H)=1$;

\ if $F^k\subset H$ then
$$
\begin{array}{ll}
F^{k-1}_{i_1}\lessdot F^k  \ \Rightarrow  &\mu(s,H)=1   \\
F^k\lessdot F^{k-1}_{i_1} \ \Rightarrow  &\\
\quad  \mu(s,H)=3 & \text{if $H$ separates $C_0$ from the first element in $s$ which is lower}  \\
\qquad \text{ than $F^k$;} &  \\
\quad \mu(s,H)=1\ &\text{otherwise.}  \\
\end{array}
$$

If $H \ \text{does not separate any two among}\ C_0,\ C,\ D
\Rightarrow \mu\leq 2$
\end{lem}

\begin{dm} The proof is very similar to that of lemma
\ref{lem:minpath}. \qed \end{dm}

\bigskip

\begin{teo}\label{thm:morsecompl2} For the local system $L(\mathcal T)$
the coefficient $\bar{u}(s)_*$ in theorem \ref{thm:morsecompl1} is
given by
$$\bar{u}(s)_*\ =\ \prod_{H\in \mathcal A}\ t_H^{m(s,H)}$$
where if $s=(F^{k-1}_{i_1},..., F^{k-1}_{i_h})$ then
$$m(s,H)\ :=\ \left[\frac{\mu(s,H)-\epsilon(C)+\epsilon(D)}{2}\right]$$
where $\epsilon(C)$ (resp. $\epsilon(D)$) holds $1$ or $0$ according
whether $H$ separetes the base chamber $C_0$ from $C$ (resp. $D$).

Therefore one always has $m(s,H)\leq 1,$ with $m(s,H)=1$\ if

i) $H\in\mathcal H(C_0,C)\cap\mathcal H(C,D)$ and
$$F^k\subset H,\ F^k\lessdot F^{k-1}_{i_1}$$
with $H$ separating $C_0$ from the first element in $s$ which is
lower than $F^k;$

ii) $H\in\mathcal H(C_0,C)\cap\mathcal H(C_,D)$ and
$$F^k\subset H,\ F^k\lessdot F^{k-1}_{i_1}.$$

In the other cases we have
$$m(s,H)\leq 1$$
if $H$ does not separates any two of the three chambers $C_0,\ C, \
D,$ otherwise
$$m(s,H)=0.$$
\end{teo}
\begin{dm} The proof follows directly from the previous lemma,
by computing, for each $s,$ the number of times the path
$\bar{u}(s)$ turns around some hyperplane. \qed
\end{dm}

\bigskip

Theorem \ref{thm:morsecompl2} gives an efficient algorithm to
compute abelian local systems in terms of the polar ordering (see
also \cite{dcohen, cohenorlik, esnschvie, kohno, libyuz, salvetti2,
schtervar, suciu, yoshinaga}).

\section{The braid arrangement}

In this section, we describe the  combinatorial gradient vector
field for the {\it braid arrangement} $\mathcal{A}\ =\ \{
H_{ij}=\{x_i=x_j\},\ 1 \leq i<j \leq n+1\}.$
Let us start with some notations.\\

\subsection{Tableaux description for the complex $\S(A_n)$}
We indicate simply by $A_n$ the symmetric group on $n+1$ elements,
acting by permutations of the coordinates. Then $\Cal{A} =
\Cal{A}(A_n)$ is the braid arrangement and $\mathbf{S}(A_n)$ is the
associated CW-complex (see \ref{S-complex}).

Given a system of coordinates in $\R^{n+1}$, we describe
$\mathbf{S}(A_n)$ through certain tableaux as follow.

Every $k$-cell $[C \prec F]$ is represented by a tableau with $n+1$
boxes and  $n+1-k$ rows (aligned on the left), filled with all the
integers in $\{1,...,n+1\}.$ There is no monotony condition on the
lengths of the rows. One has:
\medskip

\ni - $(x_1,\ldots, x_{n+1})$ is a point in $F$ iff:\\

$1.$ $i$ and $j$ belong to the same row iff $x_i=x_j$,

$2.$ $i$ belongs to a row less than the one containing $j$ iff
$x_i < x_j$;\\

\ni - the chamber $C$ belongs to the half-space $x_i < x_j$ iff:\\

$1.$ either the row which contains $i$ is less than the one
containing $j$ or

$2.$ $i$ and $j$ belong to the same row and the column which
contains $i$ is less than the one containing $j$.
\medskip

Notice that the geometrical action of $A_n$ on the stratification
induces a natural action on the complex $\S,$ which, in terms of
tableaux, is given by a left action of $A_n$: $\sigma. \ T$ is the
tableau with the same shape as $T,$ and with entries permuted
through $\sigma.$

\subsection{Construction of singular tableaux and polar ordering}

In this part we use theorem \ref{thm:algo} constructing and ordering
"singular" tableaux,  corresponding to codimensional-$k$ facets
which intersect $V_k.$ We give both an algorithmic construction,
generating bigger dimensional tableaux from the lower dimensional
ones, and an explicit one.

 Denote by $\Tbf(A_n)$ the set of
"row-standard" tableaux, i.e. with entries increasing along each
row. Each facet in $\St$ corresponds to an equivalence class of
tableaux, where the equivalence is up to row preserving
permutations. So there is a $1-1$ correspondence between $\Tbf(A_n)$
and the set of facets in $\Cal{A}(A_n)$. Let
$\Tbf^{\mathbf{k}}(A_n)$ be the set of tableaux of dimension $k$
(briefly, $k-\text{tableaux}$), i.e. tableaux with exactly $n+1-k$
rows.  Moreover, write $T \prec T'$ iff $F \prec F'$, where the
tableaux $T$ and $T'$ correspond respectively to $F$ and $F'$. Our
aim is to give a polar ordering on $\Tbf(A_n).$

\begin{df}[Moving Function] Fixed an integer $1 \leq r \leq n+1$, for
each $0 \leq j \leq n-k$, define the $\mathrm{moving\ function}$
$$M_{j,r}:\Tbf^{\mathbf{k+1}}(A_n) \too \Tbf(A_n),$$
where the tableau $M_{j,r}(T^{k+1})$ is obtained from $T^{k+1}$
moving the entry $r$ to the $j$-th row. Case $j=0$ means that $r$
becomes the only entry of the first row in
$M_{j,r}(T^{k+1})$.\\
\end{df}

Of course, if $r$ is the unique element of its row, moving $r$ makes
the preceding and following rows to become adjacent. So, the number
of rows of the new tableau can increase or decrease by $1,$ or it
can remain equal (when the row of $r$ has at least two elements and
$j>0$).  Given a tableau $T^k,$ where $r$ is in the $i$-th row, we
define the set of tableaux $M_r(T^k)=\{M_{j,r}(T^k)\}_{0 < j < i}$.
We assign to $M_r(T^k)$ the reverse order with respect to $j$.

Let us consider the natural projection $p_{n,m}:\Tbf(A_n) \too
\Tbf(A_m)$ obtained by forgetting the entries $r \geq m+2$ in each
tableau ("empty" rows are deleted). For any $T \in \Tbf(A_n)$ denote
by $m_T$ the minimum integer $1 \leq m \leq n$ such that
$p_{n,m}(T)$ preserves the dimension of $T$. So, each $j> m_T+1$\ is
the unique element of its row.

\begin{df}[T-Blocks] Let $T$ be a $k$-tableau in
$\Tbf(A_n)$ and $e_i(T)$ the first entry of its $i$-th row. Then,
if\ $m_T<n,$  for any integer $m_{T}+1-k<h\leq n+1-k$ we define a
new ordered set
\begin{equation}\label{t-blocco}
\Cal{Q}_{n,h}(T)= \bigcup_{m_{T}+1-k < i \leq h}M_{e_i(T)}(T_{i-1}),
\hspace{0,2cm} T_{m_{T}-k+1}=T \mbox{ and }
T_i=M_{0,e_i(T)}(T_{i-1}),
\end{equation}
where $M_{e_i(T)}(T_{i-1})$ are already ordered and tableaux in
$M_{e_i(T)}(T_{i-1})$ are less than tableaux in
$M_{e_j(T)}(T_{j-1})$ iff $i < j$.

\end{df}

Let $T \in \Tbf(A_n)$ be a tableau representing a facet $F$. The
symmetry in $\R^{n+1}$ with respect to the affine subspace generated
by $F$ preserves the arrangement, so it induces an involution $r_T$
on $\Tbf(A_n).$ Given a $k$-tableau $T \in \Tbf^{\mathbf{k}}(A_n)$
with $m_T<n,$ let $\Cal{Q}_{n,h}(T)\ =\ \{T_i\}_{1 \leq i \leq p},$
where the indices follow the ordering introduced in the previous
definition.  If we consider a $k$-tableau $\overline{T}$ then we can
define recursively $\overline{T}_i$ as follow:
\begin{enumerate}
\item $\overline{T}_{1}=\overline{T}$
\item $\overline{T}_i=r_{T_i}\overline{T}_{i-1}$ if
$T_i \succ \overline{T}_{i-1}$, $\overline{T}_i=\overline{T}_{i-1}$
otherwise.
\end{enumerate}
Denote the last tableau $\overline{T}_p$  by $r_{\Cal{Q}_{n,h}(T)}(\overline{T}).$\\

Let $i_{m,n}:\Tbf(A_{m}) \too \Tbf(A_n)$ be the natural inclusion
map, i.e. $i_{m,n}(T)$ is obtained by attaching to $T$ exactly $n-m$
rows of lenght one having entries $m+2,\ ...\ ,n+1$ increasing along
the first column.

Let $\pi_0(A_n)$ be the set given by the identity $0$-tableau (i.e.,
one column with growing entries); we define $\pi_{k+1}(A_n) \subset
\Tbf^{\mathbf{k+1}}(A_n)$ as the image of the map:
\begin{equation}\label{kspaziocritico}
\begin{split}
\Cal{\overline{Q}}_{n,n+1-k}: &\pi_{k}(A_{n-1}) \too
\Tbf^{\mathbf{k+1}}(A_{n}),\\ &T_i \too \Cal{Q}_{n,n+1-k}(T_{i,i})
\end{split}
\end{equation}
where $T_{i,1}=i_{n-1,n}(T_i)$ and
$T_{i,j}=r_{\Cal{Q}_{n,n+1-k}(T_{j-1,j-1})}(T_{i,j-1})$ for $j \leq
i$.

We inductively order  $\pi_k(A_{n})$ by requiring that the map in
\ref{kspaziocritico} is order preserving  and using the ordering of
the T-blocks involved.

\begin{rmk}\label{rifletto} Remark that the $k$-tableau $T_{i,i}$
in the above definition is $T_{i,i}=r_{T}(T_{i,1})$ where
$T=i_{m_{T_i},n}(T^{m_{T_i}})$ and $T^{m_{T_i}}$ is the unique
tableau (having only one row) of
$\Tbf^{\mathbf{m_{T_i}}}(A_{m_{T_i}})$.
\end{rmk}

Now let us describe directly tableaux $T^k$ in $\pi_k(A_n)$. Define
the following operations between tableaux:
\begin{enumerate}
\item $T * T^{\prime}$ is the new tableau obtained by attaching vertically
$T^{\prime}$ below $T.$
\item $T *_i h$ is the tableau obtained by
attaching the one-box tableau with entry $h$ to the $i$-th row of
$T$.
\item $T^{op}$ is the tableau obtained from $T$ by reversing the row order.
Notice that $(T * T^{\prime})^{op}=T^{\prime op} * T^{op}$.
\end{enumerate}
Let us fix $k$ integers $1 < j_1 < \cdots < j_k \leq n+1$ and, for
$1 \leq h \leq k+1$, let $T_h$ be the $0$-tableau (= one-column
tableau) with entries $J_h=\{j_{h-1}+1, \ldots , j_h-1\}$ in the
natural order (set $j_0=0,\ j_{k+1}=k+2$).\\
Then, for any suitable choice of integers $i_1,\ldots ,i_k$ we
define a $k$-tableau:
\begin{equation}\label{tabdipi}
T^k=((\cdots ((((T_1^{op}*_{i_1} j_1) * T_2)^{op}*_{i_2} j_2) *
T_3)^{op} \cdots)^{op}*_{i_k} j_k) * T_{k+1}.
\end{equation}

\begin{prop}\label{costruzione}A $k$-tableau in $\Tbf(A_n)$ is
in $\pi_k(A_n)$ iff it is of the form (\ref{tabdipi}). Moreover, the
order in $\pi_k(A_n)$ is the one induced by lexicographic order
between sequences of pairs $((j_1,i_1),\ldots ,(j_k,i_k))$, where
$(j_t,i_t) < (j_t^{\prime},i_t^{\prime})$ iff either $j_t <
j_t^{\prime}$ or $j_t = j_t^{\prime}$ and $i_t > i_t^{\prime}$.
\end{prop}

\textbf{Proof.} The proof is by induction on the dimension $n$ of
$\Cal{A}(A_n)$.\\
The result holds trivially for $n=1$.\\
Let $T^k$ be a tableau in $\Tbf(A_n)$ such that the $(n+1-k)$-th row
has length one and entry $n+1$. Then, by construction, $T^k \in
\pi_k(A_n)$ iff $p_{n,n-1}(T^k) \in \pi_k(A_{n-1})$ and proof comes
by inductive hypothesis.\\
Otherwise $j_k=n+1$, i.e.
$$ T^k=((\cdots (((T_1^{op}*_{i_1} j_1) * T_2)^{op}*_{i_2} j_2 * T_3)^{op}
\cdots)^{op} *_{i_{k-1}} j_{k-1} * T_k)^{op}*_{i_k} (n+1).$$ If we
define a $(k-1)$-tableau as
$$T^{k-1}=(\cdots (((T_1^{op}*_{i_1} j_1) * T_2)^{op}*_{i_2} j_2
* T_3)^{op} \cdots)^{op} *_{i_{k-1}} j_{k-1} * T_k $$ then
$T^{k-1} \in \pi_{k-1}(A_{n-1})$ by induction and $T^k\in
\Cal{\overline{Q}}_{n,n+1-(k-1)}(T^{k-1})$ by construction, i.e.
$T^k \in \pi_k(A_n)$. Since, by construction, given $T, T^{\prime}$
belonging to $\Cal{\overline{Q}}_{n,n+1-(k-1)}(T^{k-1}),$ one has
that $T$ is lower than  $T^{\prime}$ iff either $j_k < j_k^{\prime}$
or $j_k = j_k^{\prime}$ and $i_k > i_k^{\prime}$ then the proof
arises from inductive hypothesis. \qed

\bigskip

Let us consider the subset $\Cal{U}^k(A_n)$ of rank-$k$ elements in
the lattice $L(\Cal{A}(A_n))$ (see \cite{orlik_terao}): in other
words, the set of codimensional-$k$ intersections of hyperplanes
from $\Cal A$. The support of the facet represented by $T^k\in
\Tbf^{\mathbf{k}}(A_n)$ is denoted by $\mid T^k \mid\in
\Cal{U}^k(A_n).$

By arguments similar to those used in the proof of proposition
\ref{costruzione}, one obtains the following result.

\medskip

\begin{lem}\label{spazispasi} $\pi_k(A_n)$ is a complete system of
representatives for $\Cal{U}^k(A_n)$, i.e. any affine space in
$\Cal{U}^k(A_n)$ is the support of a tableau in $\pi_k(A_n)$ and any
two $k$-tableaux in $\pi_k(A_n)$ have different supports. \qed
\end{lem}

\medskip

\begin{rmk} It follows that the cardinality of $\pi_k(A_n)$ is
the number of $k$-codimensional subspaces of the intersection
lattice $L(\Cal{A}(A_n))$, i.e. the Stirling number $S(n+1,n+1-k)$
(see \cite{orlik_terao}).
\end{rmk}

Now let us prove that tableaux in $\pi_k(A_n)$ describe critical
cells of $\S^k(A_n)$ with respect to a suitable system of polar
coordinates.

\begin{prop}\label{k-sezione} There exists a system of polar
coordinates, generic with respect to $\mathcal{A}(A_n),$ such that a
codimensional-$k$ facet $F$ meets the $V_k$ space iff the tableau
which represent $F$ is in $\pi_k(A_n)$. Moreover, the induced order
between codimensional-$k$ facets intersecting $V_k$ equals that
introduced before for $\pi_k(A_n)$.
\end{prop}

\textbf{Proof.} We start defining
$\Cal{A}(A_{n-1}^{n})=i_{n-1,n}(\Cal{A}(A_{n-1}))$,
$\Cal{A}(A_{n-1}^{n})^{c}=\Cal{A}(A_n) \setminus
\Cal{A}(A_{n-1}^{n}).$ Let also $\pi_{k-1}(A_{n-1}^n)=
i_{n,n-1}(\pi_{k-1}(A_{n-1})).$

The proof is by double induction on the dimension $n$ of
$\Cal{A}(A_n)$ and the dimension $k$ of sections $V_k$.  The result
holds trivially for $n=1,\ 2$ and also for $k=0$ and any $n.$

By induction, it is possible to find a system of generic polar
coordinates $V'_0,...,V'_{n}$ in $\R^n$  which verifies the theorem
for $A_{n-1}.$ By using arguments similar to those used in section
\ref{sec:polarcoord} one can embed this system to a generic one
$V_0,...,V_n, \ V_{n+1}=\R^{n+1}$ for $A_{n},$ where the embedding
is compatible with $i_{n,n-1}$ (i.e., it takes $\Cal A(A_{n-1})$
inside $\Cal A(A^n_{n-1})$).

By induction on $k$, we assume that the system verifies the
assertion until codimensional-$(k-1)$ facets.

Let $\Cal{L}_k(\pi_{k-1}(A_{n-1}^n))$ be the set of all affine lines
realized as intersections between $V_k$ and $\Cal{U}^{k-1}(\Cal
A(A_{n-1}^n))$. By lemma \ref{spazispasi} any line $L_i$ in
$\Cal{L}_k(\pi_{k-1}(A_{n-1}^n))$ lies in the support of one and
only one tableau $T_i^{k-1} \in \pi_{k-1}(A_{n-1}^{n}).$

Now notice that for any $T_i^{k-1} \in \pi_{k-1}(A_{n-1}^{n})$, the
last row is composed only by the entry $(n+1).$ Moreover, by remark
\ref{rifletto}, the tableau $T_{i,i}^{k-1}$ is obtained from
$T_i^{k-1}$ without moving the entry $n+1$. Then, by construction,
$\pi_k(A_{n-1}^{n})$ is given by the ordered union of
$\Cal{{Q}}_{n,n-(k-1)}(T_i^{k-1})$ for $T_i^{k-1} \in
\pi_{k-1}(A_{n-1}^{n})$.

Therefore (by induction) the line $L_i$  intersects in order all
$k$-facets represented by $\Cal{{Q}}_{n,n-(k-1)}(T_i^{k-1})$ and,
after that, all hyperplanes in $\Cal{A}(A_{n-1}^{n})^c$. Obviously
these last intersections have to be along a gallery of $k$-tableaux
starting from the $(k-1)$-tableau
$$\tilde{T}^{k-1}_i:=r_{\Cal{Q}_{n,n-(k-1)}(T_i^{k-1})}(T_i^{k-1}).$$ But
$M_{e_{(n+1)-(k-1)}(\tilde{T}^{k-1}_i)}(\tilde{T}^{k-1}_i)=M_{n+1}(\tilde{T}^{k-1}_i)$
is the only choice in order to have a gallery throughout hyperplanes
in $\Cal{A}(A_{n-1}^{n})^c$ and starting from $\tilde{T}^{k-1}_i$.

This proves the first statement of the proposition.

According to definition \ref{df:minpoint} let
$$P_{F^k_{i,h}}\ :=\ clos(F^k_{i,h})\cap V_k$$
where $F^k_{i,h}$ is the facet represented by the tableau $T_{i,h}^k
\in\Cal{{Q}}_{n,(n+1)-(k-1)}(T_i^{k-1}).$ By definition
\ref{df:maindef} we need to understand the ordering of such points
$P's$.

 By the above considerations and the inductive hypothesis it follows that in
$\pi_k(A_{n-1}^{n})$ one has
$$P_{F^k_{i_1,h_1}} \lessdot P_{F^k_{i_2,h_2}}$$
iff the pair $(i_2,h_2)$ follows the pair $(i_1,h_1)$ according to
the lexicographic ordering.  By simple geometric considerations this
lexicographic ordering is preserved when we pass to $\pi_k(A_{n}).$
But this corresponds exactly to the ordering which we defined before
for $\pi_k(A_{n}).$  \qed

\begin{rmk}\label{osservaz} By theorem
\ref{thm:algo}, we can reconstruct the ordering of $\Tbf(A_n)$ from
that of $\pi_k(A_n),\ k=0,...n.$

In order to identify critical cells of $\S(\Cal{A}(A_n))$ we just
apply theorem \ref{T6}.

\end{rmk}

\providecommand{\bysame}{\leavevmode\hbox
to3em{\hrulefill}\thinspace}
\providecommand{\MR}{\relax\ifhmode\unskip\space\fi MR }
\providecommand{\MRhref}[2]{%
  \href{http://www.ams.org/mathscinet-getitem?mr=#1}{#2}
} \providecommand{\href}[2]{#2}


\begin{thebibliography}{DPSS99}







\bibitem[BZ92]{bjorner_ziegler}
A. Bjorner, G. Ziegler, \emph{Combinatorial stratifications of
complex arrangements}, Jour. Amer. Math. Soc.
  \textbf{5} (1992), 105--149.


\bibitem[Bou68]{bourbaki}
N.~Bourbaki, \emph{Groupes et algebr\`es de {L}ie}, vol. Chapters
{IV-VI},
  Hermann, 1968.


\bibitem[Co93] {dcohen}
D. Cohen, \emph {Cohomology and intersection cohomology of complex
hyperplane arrangements}, Adv. Math.  \textbf{97} (1993), no.2,
231--266

\bibitem[CO00] {cohenorlik}
D. Cohen, P. Orlik, \emph {Arrangements and local systems}, Math.
Res. Lett. \textbf{7} (2000), no.2-2, 299--316



\bibitem[DP03] {dimcapapa}
A. Dimca, S. Papadima, \emph {Hypersurface complements, Milnor
fibers and higher homotopy groups of arrangements}, Ann. of Math.
(2) \textbf{158} (2003), no.2, 473--507


\bibitem[ESV92] {esnschvie}
H. Esnault, V. Schechtman, E. Viehweg, \emph {Cohomology of local
systems on the complement of hyperplanes}, Inv. Math.  \textbf{109}
(1992), no.3, 557--661




\bibitem[Fo98] {forman}
R. Forman, \emph {Morse Theory for Cell Complexes}, Adv. in Math.
\textbf{134} (1998),  no.1,  90--145


\bibitem[Fo02]{forman1}
R. Forman, \emph{A User's guide to discrete Morse Theory}, Sem.
Lotharingien de Combinatoire \textbf{48} (2002)


\bibitem[GR89]{gelfand_rybnikov}
I.M. Gelfand, G.L. Rybnikov, \emph{Algebraic and topological
invariants of oriented matroids}, Dokl.
  \textbf{307} (1989), 791--795.


\bibitem[Ko86] {kohno}
T. Kohno, \emph {Homology of a local system on the complement of
hyperplanes}, Proc. Japan Acad.  \textbf{62}, Ser. A  (1986),
144--147


\bibitem[LY00] {libyuz}
A. Libgober, S. Yuzvinsky, \emph {Cohomology of the Orlik-Solomon
algebras and local systems }, Compositio Math. \textbf{21} (2000),
337--361


\bibitem[OT92]{orlik_terao}
P. Orlik, M. Terao,  \emph{Arrangements of hyperplanes},
  Springer-Verlag \textbf{300} (1992)



\bibitem[Ra02]{randell}
R. Randell, \emph{Morse theory, Milnor fibers and minimality of a
complex hyperplane arrangement}, Proc. Amer. Math. Soc. \textbf{130}
(2002), no. 9, 2737--2743


\bibitem[Sal87]{salvetti}
M.~Salvetti, \emph{Topology of the complement of real hyperplanes in
$\C^n$}, Inv. Math.,
  \textbf{88} (1987), no.3, 603--618.


\bibitem[Sal94]{salvetti2}
M.~Salvetti, \emph{The homotopy type of {A}rtin groups}, Math. Res.
Lett.
  \textbf{1} (1994), 567--577.


\bibitem[SS07]{salvettisette}
M.~Salvetti, S, Settepanella \emph{Combinatorial Morse theory and
minimality of hyperplane arrangements}, preprint Dipartimento di
Matematica "L. Tonelli", n. 1.328.1655, January 2007.



\bibitem[STV95] {schtervar}
V. Schechtman, H. Terao, A. Varchenko, \emph {Cohomology of local
systems and the Kac-Kazhdan condition for singular vectors}, J. Pure
Appl. Alg. \textbf{100} (1995), 93--102



\bibitem[Su02] {suciu}
A. Suciu, \emph {Translated tori in the characteristic varieties of
complex hyperplane arrangements}, Topol. and Appl.  \textbf{118}
(2002), no.1-2, 209-223



\bibitem[Za75]{zaslavsky}
T.~Zaslavsky, \emph{Facing up to arrangements: face count formulas
for partitions of space by hyperplanes}, Memoirs Amer. Math. Soc.
  \textbf{1} (1975), no.154.



\bibitem[Yo05]{yoshinaga}
M. Yoshinaga, \emph{Hyperplane arrangements and Lefschetz's
hyperplane section theorem} (2005)



\end{thebibliography}
\end{document}